\documentclass[12pt,twoside]{amsart}
\usepackage{amssymb}
\usepackage{amscd}
\usepackage[abbrev,alphabetic]{amsrefs}
\usepackage{hyperref}
\usepackage{here}
\title[On p-power freeness in positive characteristic]
{On p-power freeness in positive characteristic} 
\author{Hiromu Tanaka} 
\address{Graduate School of Mathematical Sciences, 
The University of Tokyo, 
3-8-1 Komaba, Meguro-ku, Tokyo 153-8914, JAPAN} 
\email{tanaka@ms.u-tokyo.ac.jp}
\newcommand{\red}[0]{{\operatorname{red}}}

\newcommand{\Spec}[0]{{\operatorname{Spec}}}

\newcommand{\Supp}[0]{{\operatorname{Supp}}}

\newcommand{\Ex}[0]{{\operatorname{Ex}}}

\newtheorem{thm}{Theorem}[section]
\newtheorem{lem}[thm]{Lemma}

\newtheorem{prop}[thm]{Proposition}

\newtheorem{step}{Step}

\theoremstyle{definition}

\newtheorem{conj}[thm]{Conjecture}
\newtheorem{dfn}[thm]{Definition}

\newtheorem{rem}[thm]{Remark}

\makeatletter
  
  \@addtoreset{equation}{thm}
  \makeatother

\newcommand{\MO}{\mathcal{O}}
\newcommand{\R}{\mathbb{R}}
\newcommand{\Q}{\mathbb{Q}}
\newcommand{\Z}{\mathbb{Z}}
\newcommand{\F}{\mathbb{F}}

\newcommand{\m}{\mathfrak{m}}

\newcommand{\mcL}{\mathcal{L}}
\begin{document}

\maketitle

\begin{abstract}
In this note, we study base point freeness up to taking $p$-power, 
which we will call $p$-power freeness. 
We first establish some criteria for $p$-power freeness as analogues of criteria for semi-ampleness. 
We then apply these results to three-dimensional birational geometry.  
\end{abstract}

\tableofcontents

\section{Introduction}

It is a fundamental problem  
to study torsion line bundles in algebraic geometry. 
%Especially, global features of algebraic varieties are reflected to the existence or non-existence of non-trivial torsion line bundles. 
If $L$ is an $m$-torsion line bundle on an algebraic variety $X$ of characteristic zero, 
then $L$ induces an \'etale cover of degree $m$. 
The same statement holds in characteristic $p>0$ when $m$ is not divisible by $p$. 
However, if $m$ is divisible by $p$, 
then the resulting cyclic cover $Y \to X$ is no longer \'etale. 
Therefore $p^e$-torsion line bundles have different feature from $\ell$-torsion line bundles for $\ell \in \Z_{>0} \setminus p\Z$.

%On the other hand, we naturally encounter $p^e$-torsion line bundles in algebraic geometry in characteristic $p>0$. 
To observe another phenomenon on $p^e$-torsion line bundles,
let us recall the Lefschetz theorem for the local Picard groups (cf. \cite[Theorem 0.1]{BdJ14}). 
Let $(A, \m)$ be an excellent normal local ring containing a field. 
For $f \in \m \setminus \{0\}$, $V:=\Spec\,A \setminus \{\m\}$, and $V_0 :=\Spec\,(A/f) \setminus \{\m\}$, 
the restriction map $\rho:{\rm Pic}(V) \to {\rm Pic}(V_0)$ satisfies the following properties. 
\begin{enumerate}
\item[(0)] Assume that $A$ is of characteristic zero. 
If ${\rm depth}_{\m}(A/f) \geq 2$, then $\rho$ is injective. 
\item[(p)] 
Assume that $A$ is of positive characteristic. 
Then $\rho$ is injective up to taking $p$-power, i.e. 
if $\rho(L) \simeq \MO_{V_0}$, then $L^{\otimes p^e} \simeq \MO_V$ for some $e \in \Z_{>0}$. 
\end{enumerate}
In this theorem, we need to care $p^e$-torsion line bundles for the case of positive characteristic, 
whilst there does not appear such non-trivial torsion line bundles 
in the corresponding statement of characteristic zero. 
One of the purposes of this note is to observe 
a similar phenomenon that appears in birational geometry.

%One of geometric meaning of torsion line bundles is reflected 
%to its associated finite covers. 
%For example, 
%if $L$ is an $m$-torsion line bundle on a variety $X$ of characteristic zero, 
%then $L$ induces an \'etale cover of degree $m$. 
%However, the same statement is no longer true in characteristic $p>0$. 
%If $L$ is an $m$-torsion line bundle for $m \in \Z_{>0} \setminus p\Z$, then it still induces an \'etale cover of degree $m$. On the other hand, if $L$ is a $p$-torsion line bundle, then the associated cover $Y \to X$ is purely inseparable. This cover is still useful, however study of such cover $Y$ is more difficult, because $Y$ might not be smooth even if so is $X$. We naturally encounter $p^e$-torsion line bundles on $X$. 

%find criteria for base point freeness. For instance, the Kawamata--Shokurov base point free theorem is one of essential results in birational geometry. On the other hand, sometimes the problem becomes easier in positive characteristic. The typical example is Keel's theorem. However, the proof of Keel's theorem uses Frobenius morphism to extending sections, hence it compels us to take $p^e$-power of given line bundles to conclude base point freeness. Actually, it is sometimes necessary to take such $p$-powers in characterisitc $p>0$ (cf. ). 

To this end, we first study the behaviour of $p^e$-torsion line bundles in general settings. 
For flexibility, 
we study a wider notion: 
base point freeness up to $p$-power, which we will call $p$-power freeness. 

\begin{dfn}\label{intro-pp-free}
Let $f:X \to Y$ be a proper morphism of noetherian $\F_p$-schemes. 
Let $L$ be an invertible sheaf on $X$. 
\begin{enumerate}
\item 
We say that $L$ is {\em $f$-free} 
if the induced homomorphism
\[
f^*f_*L \to L
\]
is surjective. 
If $Y=\Spec\,k$ for a field $k$, then we simply say that $L$ is {\em free}.
\item 
We say that $L$ is {\em p-power $f$-free} 
or {\em p-power free over $Y$} if 
there exists a positive integer $e \in \Z_{>0}$ such that 
$L^{\otimes p^e}$ is $f$-free. 
If $Y=\Spec\,k$ for a field $k$, then we simply say that $L$ is 
{\em p-power free}.
\end{enumerate}
\end{dfn}

There are some criteria for semi-ampleness 
which hold only in positive characteristic. 
A typical result is Keel's theorem \cite[Theorem 0.2]{Kee99}. 
Another example is the equivalence between 
the relative semi-ampleness and the fibrewise semi-ampleness \cite[Theorem 1.1]{CT}. 
It is remarkable that these two semi-ampleness criteria 
have analogous statements for $p$-power freeness as follows. 

\begin{thm}[Theorem \ref{t-rel-pp-free1}]\label{intro-rel-pp-free1}
Let $f:X \to Y$ be a projective morphism of noetherian $\F_p$-schemes. 
Let $L$ be an invertible sheaf on $X$. 
If $L|_{X_y}$ is $p$-power free for any point $y \in Y$, 
then $L$ is $p$-power $f$-free. 
\end{thm}

\begin{thm}[Theorem \ref{t-Keel}]\label{intro-Keel}
Let $f:X \to Y$ be a projective morphism of noetherian $\F_p$-schemes. 
Let $L$ be an $f$-nef invertible sheaf on $X$ and 
let $g:\mathbb E_f(L) \hookrightarrow X \xrightarrow{f} Y$ be the induced morphism. 
Then $L$ is $p$-power $f$-free if and only if $L|_{\mathbb E_f(L)}$ is $p$-power $g$-free. 
\end{thm}

We then apply these criteria to birational geometry in positive characteristic. 
Let us first recall the Kawamata--Shokurov base point free theorem in characteristic zero \cite[Theorem 3-1-1]{KMM87}. 

\begin{thm}[Kawamata--Shokurov]\label{intro-KS-0}
Let $k$ be a field of characteristic zero. 
Let $(X, \Delta)$ be a klt pair over $k$ and 
let $f:X \to Z$ be a projective $k$-morphism to a quasi-projective $k$-scheme $Z$. 
Let $L$ be an $f$-nef Cartier divisor such that $L-(K_X+\Delta)$ is $f$-nef and $f$-big. 
Then there exists a positive integer $m_0$ such that 
$mL$ is $f$-free for any integer $m \geq m_0$. 
\end{thm}

It is known that the same statement is no longer true in positive characteristic \cite[Theorem 1.2]{Tan1}. 
More specifically, 
over an algebraically closed field $k$ of characteristic $p \in \{2, 3\}$, 
there exist a three-dimensional klt pair $(X, \Delta)$, 
a projective morphism $f:X \to Z$ to a smooth curve $Z$, and 
an $f$-numerically trivial Cartier divisor $L$ on $X$ such that $L-(K_X+\Delta)$ is $f$-ample and $L\not\sim_f 0$. 
On the other hand, it holds that $pL \sim_f 0$ for this example. 
Then it is tempting to hope that $L$ is $p$-power free in the case of positive characteristic.

\begin{conj}\label{c-KS-bpf-p}
Let $k$ be a field of characteristic $p>0$. 
Let $(X, \Delta)$ be a klt pair over $k$ and 
let $f:X \to Z$ be a projective $k$-morphism to a quasi-projective $k$-scheme $Z$. 
Let $L$ be an $f$-nef Cartier divisor on $X$ such that $L-(K_X+\Delta)$ is $f$-nef and $f$-big. 
Then $L$ is $p$-power $f$-free, i.e. 
there exists a positive integer $e$ such that $p^eL$ is $f$-free. 
\end{conj}

\begin{rem}
Assume that $\dim X=2$. 
Then the same statement of Theorem \ref{intro-KS-0} holds 
for the case when $k$ is a perfect field of characteristic $p>0$ 
(Lemma \ref{l-pp-2dim}). 
However, if $k$ is allowed to be an imperfect field, 
then the same statement of Theorem \ref{intro-KS-0} no longer holds \cite[Theorem 1.4]{Tan1}. 
On the other hand, Conjecture \ref{c-KS-bpf-p} is known to hold when $X$ is a surface over an imperfect field (cf. Remark \ref{r-pp-2dim}). 
\end{rem}

In this note, we establish the following partial solution 
(Theorem \ref{intro-bpf-pp}) 
as an application of our criteria for $p$-power freeness. 
%(Theorem \ref{intro-rel-pp-free1} and Theorem \ref{intro-Keel}). 
Note that our proof depends on the recent development of minimal model program for threefolds of characteristic $p>5$ 
(\cite{HX15}, \cite{CTX15}, \cite{Bir16}, \cite{BW17}, \cite{Wal18}, \cite{HNT}). 

\begin{thm}[Theorem \ref{t-bpf-pp}]\label{intro-bpf-pp}
Let $k$ be a perfect field of characteristic $p>5$. 
Let $(X, \Delta)$ be a three-dimensional klt pair and 
let $f:X \to Z$ be a projective surjective $k$-morphism to a quasi-projective $k$-scheme $Z$. 
Let $L$ be an $f$-nef Cartier divisor on $X$ such that $L-(K_X+\Delta)$ is $f$-nef and $f$-big. 
Assume that either 
\begin{enumerate}
\item $\dim Z \geq 1$, or 
\item $\dim Z=0$ and $L \not\equiv 0$. 
\end{enumerate}
Then $L$ is $p$-power $f$-free.  
\end{thm}

%\begin{rem}\label{r-KS-bpf-p}
%Our proof of Theorem \ref{intro-bpf-pp} ($=$Theorem \ref{t-bpf-pp}) 
%shows that Conjecture \ref{c-KS-bpf-p} is a consequence of the following conjectures. 
%\begin{enumerate}
%\item Minimal model program for log canonical pairs. 
%\item Existence of log resolutions. 
%\item The normalisation of the coefficient one part of a plt pair is a universal homeomorphism. 
%\item Conjecture \ref{c-KS-bpf-p} for the case when $Z=\Spec\,k$ and $L\equiv 0$. 
%\end{enumerate}
%Note that if the base field is perfect, then (1) implies (3) in dimension three \cite[Theorem 3.14]{GNT}. 
%\end{rem}

\begin{rem}
Under the same assumptions as the ones of Theorem \ref{intro-bpf-pp}, 
Bernasconi proves that if 
$p \gg 0$, $\dim Z \geq 1$, and $L$ is a Cartier divisor on $X$ such that $L \equiv_f 0$, 
then $L \sim_f 0$ \cite{Ber}. 
%Although his assumption is stronger, also the conclusion is 
\end{rem}

\subsection{Overviews of proofs}

We now overview how to prove the main theorems. Let us first discuss the criteria for $p$-power freeness (Theorem \ref{intro-rel-pp-free1}, Theorem \ref{intro-Keel}). 
Theorem \ref{intro-Keel} follows from the same argument as the proof of Keel's theorem given  by \cite{CMM14}. 
Theorem \ref{intro-rel-pp-free1} is the $p$-power free version of \cite[Theorem 1.1]{CT}. 
Although the same argument as in \cite[Theorem 1.1]{CT} does not work for our case, Theorem \ref{intro-rel-pp-free1} follows from a combination of \cite[Theorem 1.1]{CT} and a recent result by Bhatt--Scholze \cite{BS17}.

%\footnote{Birational decent should be mentioned here}

We now give an overview of how to show Theorem \ref{intro-bpf-pp}. 
The first step is to establish a birational version as follows.

\begin{thm}[Theorem \ref{t-birat-torsion}]\label{intro-t-birat-torsion}
Let $k$ be a perfect field of characteristic $p>5$. 
Let $f:X \to Y$ be a projective birational $k$-morphism of quasi-projective normal threefolds over $k$. 
Assume that there exists an effective $\R$-divisor $\Delta_Y$ on $Y$ such that 
$(Y, \Delta_Y)$ is klt. 
Let $L$ be a Cartier divisor on $X$ such that $L \equiv_f 0$. 
Then there exist a positive integer $e$ and a Cartier divisor $L_Y$ on $Y$ such that $L \sim p^eL_Y$. 
\end{thm}

%and Theorem \ref{intro-bpf-pp}. 
Roughly speaking, the proof of Theorem \ref{intro-t-birat-torsion} is 
done by resolution of singularities and minimal model program. 
Thus, a crucial part is the case when $f$ is an extremal contraction of pl-type (cf. Step \ref{s1-birat-torsion} of Theorem \ref{t-birat-torsion}). 
In this case, the problem is reduced to the case of dimension two 
by using Theorem \ref{intro-Keel}. 
%The other parts of the proof of Theorem \ref{t-birat-torsion} are settled 

In order to prove Theorem \ref{intro-bpf-pp}, we utilise a variant of Theorem \ref{intro-rel-pp-free1} (Theorem \ref{t-rel-pp-free2}), 
so that it is enough to show that $L|_{X_z}$ is $p$-power free for any closed point $z \in Z$. 
Theorem \ref{intro-t-birat-torsion} enables us to take birational model changes. 
Then we may assume that $T:=(X_z)_{\red}$ is a surface of Fano type, 
i.e. there exists an effective $\Q$-divisor $\Delta_T$ such that $(T, \Delta_T)$ 
is klt and $-(K_T+\Delta_T)$ is ample 
(cf. Step \ref{s1-bpf-pp} of Theorem \ref{t-bpf-pp}). 
Then $L|_{(X_z)_{\red}}$ is trivial, hence 
also $p^eL|_{X_z}$ is trivial for some $e \in \Z_{>0}$. 
In particular, $L|_{X_z}$ is $p$-power free. 
For more details, see Section \ref{s-MMP}.

%Since both the proofs are quite similar, let us give an outline of the proof of Theorem \ref{intro-bpf-pp} using Theorem \ref{intro-bpf-pp}. 
%The first step is to handle the case when $f$ is an extremal contraction of pl-type (cf. Step %\ref{s1-birat-torsion} of Theorem \ref{t-birat-torsion}). 
%In this case, the problem is reduced to the case of dimension two by using Theorem \ref{intro-Keel}. 

\textbf{Acknowledgements:} 
The author would like to thank Fabio Bernasconi 
for useful comments. 
The author was funded by the Grant-in-Aid for Scientific Research (KAKENHI No. 18K13386).

\section{Preliminaries}

\subsection{Notation}\label{ss1-prelim}

\begin{enumerate}
\item 
We will freely use the notation and terminology in \cite{Har77} 
and \cite{Kol13}. 
\item 
For a scheme $X$, its {\em reduced structure} $X_{\red}$ 
is the reduced closed subscheme of $X$ such that the induced morphism 
$X_{\red} \to X$ is surjective. 
\item 
For a field $k$, 
we say that $X$ is a {\em variety over} $k$ or a $k$-{\em variety} if 
$X$ is an integral scheme that is separated and of finite type over $k$. 
We say that $X$ is a {\em curve} over $k$ or a $k$-{\em curve} 
(resp. a {\em surface} over $k$ or a $k$-{\em surface}, 
resp. a {\em threefold} over $k$) 
if $X$ is a $k$-variety of dimension one (resp. two, resp. three).  
\item 
Let $f: X \to Y$ be a morphism of noetherian schemes. 
We say that $f$ is {\em projective} 
if there exists a closed immersion $X \hookrightarrow \mathbb P^n_Y$ over $Y$ 
for some $n \in \Z_{>0}$. 
This definition coincides with the one in \cite[page 103]{Har77}, 
but differs from the one given by Grothendieck \cite[D\'efinition 5.5.2]{Gro61}. 
On the other hand, their definitions coincide in many cases 
(cf. \cite[Section 5.5.1]{FGAex}). 
\item 
A morphism $f:X \to Y$ of schemes {\em has connected fibres} 
if $X \times_Y \Spec\,L$ is either empty or connected for any field $L$ and any morphism $\Spec\,L \to Y$. 
\end{enumerate}

\subsection{Properties of divisors}\label{ss2-prelim}

Let $f: X\to Y$ be a proper morphism of noetherian schemes 
and let $L$ be an invertible sheaf on $X$. 

\begin{enumerate}
\item 
$L$ is $f$-{\em nef} if 
for any field $K$, morphism $\Spec\,K \to Y$, and 
a curve $C$ on $X \times_Y \Spec\,K$, 
the inequality $\alpha^*L \cdot C \geq 0$ 
holds for the induced morphism $\alpha: X \times_Y \Spec\,K \to X$. 
\item 
$L$ is $f$-{\em numerically trivial} if 
both $L$ and $L^{-1}$ are $f$-nef. 
\item\label{ss2-prelim-free} 
$L$ is $f$-{\em free} if the natural homomorphism $f^*f_*L\to L$ is surjective.  
In particular, if $L$ is $f$-free then it induces a morphism $X\to \mathbb P(f_*L)$ over $Y$. 
\item 
 $L$ is {\em $f$-very ample} if it is $f$-free and the induced morphism $X\to \mathbb P(f_*L)$ is a closed immersion. 
\item  
$L$ is $f$-{\em semi-ample} (resp. {\em  $f$-ample}) if $L^{\otimes m}$ is $f$-free (resp. $f$-very ample) for some positive integer $m$. 

\item $L$ is $f$-{\em weakly big} if there exist an $f$-ample invertible sheaf $A$ on $X$ 
and a positive integer $m$ such that if $g: X_{\red} \to Y$ denotes the induced morphism, then  
\[
g_*((L^{\otimes m}\otimes_{\MO_X} A^{-1})|_{X_{\red}})\neq 0.
\]
Assume that $X$ is normal. $L$ is $f$-{\em big} if, for any connected component $X'$ of $X$, the restriction  $L|_{X'}$ is $f'$-weakly big, where $f':X' \to Y$ denotes the induced morphism. 
\item 
If $L$ is $f$-nef, the $f$-{\em exceptional locus} of $L$, denoted by $\mathbb E_f(L)$, is defined as the union of all the reduced closed subschemes  $V\subset X$ such that $L|_V$ is not $f|_V$-weakly big. 
It is known that $\mathbb E_f(L)$ is a closed subset of $X$ \cite[Lemma 2.18]{CT}. 
We consider $\mathbb E_f(L)$ as a reduced closed subscheme of $X$. 
\end{enumerate}

\section{$p$-power freeness}\label{s-pp-free}

In this section, we first introduce $p$-power freeness for invertible sheaves (Subsection \ref{ss1-pp-free}). 
In Subsection \ref{ss2-pp-free}, 
we prove that the relative $p$-power freeness is equivalent to the fibrewise $p$-power freeness (Theorem \ref{t-rel-pp-free1}). 
In Subsection \ref{ss3-pp-free}, we establish the $p$-power free version of Keel's theorem (Theorem \ref{t-Keel}).

\subsection{Definition}\label{ss1-pp-free}

\begin{dfn}\label{d-pp-free}
Let $f:X \to Y$ be a proper morphism of noetherian $\F_p$-schemes. 
Let $L$ be an invertible sheaf. 
\begin{enumerate}
\item 
We say that $L$ is {\em p-power $f$-free} 
or {\em p-power free over $Y$} if 
there exists a positive integer $e \in \Z_{>0}$ such that 
$L^{\otimes p^e}$ is $f$-free (for definition of being $f$-free, see Subsection \ref{ss2-prelim}(\ref{ss2-prelim-free})). 
If $Y=\Spec\,k$ for a field $k$, then we simply say that $L$ is 
{\em p-power free}.
\item 
The {\em $p$-power} $f$-{\em base locus} of $L$ is defined as  the following closed subset of $X$: 
\[
\mathbb B^p_f(L)=\bigcap_{e=0}^{\infty} \Supp ~ {\rm Coker} (f^*f_*L^{\otimes p^e}\to L^{\otimes p^e}).
\]
\end{enumerate}
Note that $L$ is $p$-power $f$-free if and only if $\mathbb B^p_f(L)=\emptyset$.
\end{dfn}

\begin{lem}\label{l-pp-ff-descent}
Let 
\[
\begin{CD}
X' @>\alpha>> X\\
@VVf'V @VVfV\\
Y' @>\beta>> Y
\end{CD}
\]
be a cartesian diagram consisting of morphisms of noetherian $\F_p$-schemes, 
where $f$ is proper and $\beta$ is faithfully flat.  
Let $L$ be an invertible sheaf on $X$. 
Then $L$ is $p$-power $f$-free if and only if $\alpha^*L$ is 
$p$-power $f'$-free. 
\end{lem}

\begin{proof}
The assertion follows from the corresponding statement for usual freeness. 
\end{proof}

\begin{lem}\label{l-pp-univ-homeo}
Let 
\[
f':X' \xrightarrow{\pi} X \xrightarrow{f} Y
\]
be proper morphisms of noetherian $\F_p$-schemes. 
Assume that $\pi$ is a surjective morphism which have connected fibres.  
Let $L$ be an invertible sheaf on $X$. 
Then $L$ is $p$-power $f$-free if and only if $\pi^*L$ is $p$-power $f'$-free. 
\end{lem}

\begin{proof}
Taking the Stein factorisation of $\pi$, 
we may assume that $\pi_*\MO_{X'}=\MO_X$ or $\pi$ is finite. 
If $\pi_*\MO_{X'}=\MO_X$, then the assertion follows from the corresponding statement for usual freeness. 
Thus we may assume that $\pi$ is finite. 
In this case, $\pi$ is a finite universal homeomorphism.  
Assuming that $\pi^*L$ is $p$-power $f'$-free, 
it suffices to show that $L$ is $p$-power $f$-free. 
It follows from \cite[Proposition 6.6]{Kol97} that there exists a positive integer $e$ such that the $e$-th iterated absolute Frobenius morphism $F^e:X \to X$ factors through $\pi$: 
\[
F^e:X \xrightarrow{\rho} X' \xrightarrow{\pi} X. 
\]
Since $\pi^*L$ is $p$-power free over $Y$, 
so is its pullback 
\[
L^{\otimes p^e} = (F^e)^*L = \rho^*\pi^*L. 
\]
Hence, $L$ is $p$-power free over $Y$. 
\end{proof}

\subsection{Fibrewise $p$-power freeness}\label{ss2-pp-free}

In this subsection, we prove that 
the relative $p$-power freeness is equivalent to the fibrewise $p$-power freeness (Theorem \ref{t-rel-pp-free1}). 
If the base field is uncountable, then we obtain a stronger criterion 
(Theorem \ref{t-rel-pp-free2}).

\begin{thm}\label{t-rel-pp-free1}
Let $f:X \to Y$ be a projective morphism of noetherian $\F_p$-schemes. 
Let $L$ be an invertible sheaf on $X$. 
If $L|_{X_y}$ is $p$-power free for any point $y \in Y$, 
then $L$ is $p$-power $f$-free. 
\end{thm}

\begin{proof}
%Note that $L$ is $p$-power $f$-free if and only if,  
%for any point $y \in Y$, $Y':=\Spec\,\MO_{Y, y}$, and $X':=X \times_Y Y'$, 
%the induced morphism $f':X' \to Y'$ is $p$-power $f'$-free. 
%Replacing $f:X \to Y$ by $f':X' \to Y'$ for $y \in Y$, 
%we may assume that $\dim Y<\infty$. 
%We prove the assertion by induction on $\dim Y$. 
%If $\dim Y=0$, then there is nothing to show. 
%Assume $0<\dim Y<\infty$. 
It follows from \cite[Theorem 1.1]{CT} that $L$ is $f$-semi-ample. 
Thus, there exist a positive integer $m$, projective morphisms 
\[
f:X \xrightarrow{g} Z \xrightarrow{h} Y, 
\]
and an ample invertible sheaf $A_Z$ on $Z$ such that $L^{\otimes m} \simeq g^*A_Z$ 
and $g_*\MO_X=\MO_Z$. 
Since $L|_{X_y}$ is $p$-power free for any point $y \in Y$, 
$L_{X_z}$ is $p$-power free for any point $z \in Z$. 
By noetherian induction, we can find $e \in \Z_{>0}$ such that $L^{\otimes p^e}|_{X_z}$ is free for any $z \in Z$. 
Since $L^{\otimes p^e}|_{X_z}$ is numerically trivial, 
we obtain $L^{\otimes p^e}|_{X_z} \simeq \MO_{X_z}$. 
Therefore, by \cite[Theorem 1.3]{BS17}, there exists $d \in \Z_{>0}$ 
and an invertible sheaf $L_Z$ on $Z$ such that $L^{\otimes p^d} \simeq g^*L_Z$. 
Since $L_Z$ is ample over $Y$, $L_Z$ is $p$-power free over $Y$. 
Hence, also its pullback $L^{\otimes p^d} \simeq g^*L_Z$ is $p$-power free over $Y$. 
Therefore, $L$ is $p$-power free over $Y$. 
%Replacing $f:X \to Y$ by $g:X \to Z$, 
%we may assume that $L$ is $f$-numerically trivial. 
%By noetherian induction, we can find $e \in \Z_{>0}$ such that $L^{\otimes p^e}|_{X_y}$ is free for any $y \in Y$. 
%In particular, $L^{\otimes p^e}|_{X_z} \simeq \MO_{X_z}$ for any $z \in Z$. 
\end{proof}

Although the proof of the following theorem is very similar to the one of \cite[Lemma 6.1]{CT}, 
we give a proof for the sake of completeness.  

\begin{thm}\label{t-rel-pp-free2}
Let $k$ be an uncountable field of characteristic $p>0$. 
Let $f:X \to Y$ be a projective $k$-morphism of 
schemes which are of finite type over $k$. 
Let $L$ be an invertible sheaf on $X$. 
If $L|_{X_y}$ is $p$-power free for any closed point $y \in Y$, 
then $L$ is $p$-power $f$-free. 
\end{thm}

\begin{proof}
We prove the assertion by induction on $\dim Y$.  
If $\dim Y=0$, then there is nothing to show. 
Thus, we may assume that $\dim Y>0$ and that the assertion holds if the dimension 
of the base is smaller than $\dim Y$. 
By Theorem \ref{t-rel-pp-free1}, 
it is enough to show that $L|_{X_\xi}$ is $p$-power free  
for the generic point $\xi \in Y$ of an irreducible component of $Y$. 
Replacing $Y$ by an open neighbourhood of $\xi \in Y$, 
the problem is reduced to the case when 
$Y$ is an affine irreducible scheme such that $f$ is flat. 

By the semi-continuity theorem 
\cite[Ch. III, Theorem 12.8]{Har77}, 
if $e$ is a positive integer, then there exist a positive integer $d_e$ and 
a non-empty affine open subset $U_e \subset Y$ such that 
the equation 
\[
d_e=\dim_{k(y)} H^0(X_y,L^{\otimes p^e}|_{X_y})
\]
holds for any point $y \in U_e$. Since $k$ is uncountable, 
there exists a closed point 
\[
z\in \bigcap_{m\in \Z_{>0}} U_m.
\]
As $L|_{X_z}$ is $p$-power free, 
there exists a positive integer $e_0$ such that 
$L^{\otimes p^{e_0}}|_{X_z}$ is free. 
By \cite[Ch. III, Corollary 12.9]{Har77}, 
the restriction map 
\[
H^0(f^{-1}(U_{e_0}), L^{\otimes p^{e_0}}|_{f^{-1}(U_{e_0})}) \to H^0(X_z, L^{\otimes p^{e_0}}|_{X_z})
\]
is surjective. 
Since the base locus of $L^{\otimes p^{e_0}}$ is a closed subset of $X$, it is disjoint from $X_z$. 
In particular, $L^{\otimes p^{e_0}}|_{X_{\xi}}$ is free, as desired. 
\end{proof}

\subsection{Keel's theorem}\label{ss3-pp-free}

We have the $p$-power free version of 
Keel's theorem on semi-ampleness (Theorem \ref{t-Keel}). 
For a later use, 
we also establish a variant (Proposition \ref{p-Keel-div}). 

%The purpose of this subsection is to prove the $p$-power version of Keel's theorem (Theorem \ref{t-Keel}). Although the argument is almost the same as in \cite[Subsection 2.5]{CT}, we give a proof for the reader's convenience. 

%\begin{lem}\label{l-Keel}
%Let $f; X\to Y$ be a projective surjective morphism of noetherian $\F_p$-schemes. Let $L$ be an $f$-nef invertible sheaf on $X$. Then the following hold.  
%\begin{enumerate}
%\item 
%For an $f$-ample invertible sheaf $A$, a positive integer $m$, and an element $s \in  H^0(X_{\red}, (L^{\otimes m}\otimes_{\MO_X} A^{-1})|_{X_{\red}})$, if $Z$ is the reduced closed subscheme of $X$ whose support is equal to the zero set of $s$ and  $g: Z \hookrightarrow X \xrightarrow{f} S$ denotes the induced morphism, then $\mathbb E_f(L)=\mathbb E_g(L|_Z)$.
%\item $\mathbb E_f(L)=X$ if and only if $L$ is not $f$-weakly big.
%\end{enumerate}
%\end{lem}
%\begin{proof}
%See \cite[Lemma 2.18]{CT}. 
%\end{proof}

\begin{thm}\label{t-Keel}
Let $f:X \to Y$ be a projective morphism of noetherian $\F_p$-schemes. 
Let $L$ be an $f$-nef invertible sheaf on $X$ and 
let $g:\mathbb E_f(L) \hookrightarrow X \xrightarrow{f} Y$ be the induced morphism. 
Then the equation 
\[
\mathbb B^p_f(L)=\mathbb B_g^p(L|_{\mathbb E_f(L)})
\]
holds. 
In particular, $L$ is $p$-power $f$-free if and only if $L|_{\mathbb E_f(L)}$ is $p$-power $g$-free. 
\end{thm}

\begin{proof}
We may apply the same argument as in \cite[Proposition 2.20]{CT} 
after replacing the stable base loci $\mathbb B(-)$ by 
the $p$-power base loci $\mathbb B^p(-)$ 
(for the definition of $\mathbb B^p(-)$, see Definition \ref{d-pp-free}). 
\end{proof}

\begin{prop}\label{p-Keel-div}
Let $k$ be a field of characteristic $p>0$. 
Let $f:X \to Y$ be a projective $k$-morphism 
from a normal $k$-variety $X$ to a scheme $Y$ which is of finite type over $k$. 
Let $L$ be an $f$-nef Cartier divisor. 
Assume that 
there exist an $f$-ample $\Q$-Cartier $\Q$-divisor $A$ and 
an effective $\Q$-Cartier $\Q$-divisor $E$ such that $L \equiv_f A+E$. 
If $L|_{\Supp\,E}$ is $p$-power free over $Y$, then 
$L$ is $p$-power free over $Y$. 
\end{prop}

\begin{proof}
After replacing $A$ by $A':=L-E$, we may assume that 
the equation $L=A+E$ of $\Q$-divisors holds. 
It is enough to show the inclusion $\mathbb E_f(L) \subset \Supp\,E$, 
which follows from \cite[Lemma 2.18(1)]{CT}.  
\end{proof}

\section{Application to threefolds}\label{s-MMP}

As applications of results in Section \ref{s-pp-free}, 
we prove Theorem \ref{t-birat-torsion} and Theorem \ref{t-bpf-pp}. 
We start with a result on surfaces.

\begin{lem}\label{l-pp-2dim}
Let $k$ be a perfect field of characteristic $p>0$. 
Let $(X, \Delta)$ be a two-dimensional klt pair and 
let $f:X \to Y$ be a projective $k$-morphism to a quasi-projective $k$-scheme $Y$. 
Let $L$ be an $f$-nef Cartier divisor on $X$ 
such that $L-(K_X+\Delta)$ is $f$-nef and $f$-big. 
Then there exists a positive integer $m_0$ such that 
$mL$ is $f$-free for any integer $m \geq m_0$. 
In particular, $L$ is $p$-power $f$-free. 
\end{lem}

\begin{proof}
We may assume that $k$ is algebraically closed and $f_*\MO_X=\MO_Y$. 
If either $\dim Y \geq 1$ or $\dim Y=0$ and $L \not\equiv 0$, 
then the assertion follows from \cite[Theorem 4.2]{Tan18}. 
Thus, we may assume that $\dim Y=0$ and $L \equiv 0$. 
In this case, the assertion follows from \cite[Corollary 3.6]{Tan15}. 
\end{proof}

\begin{rem}\label{r-pp-2dim}
If we allow $k$ to be an imperfect field, 
then the same statement as in Lemma \ref{l-pp-2dim} 
no longer holds \cite[Theorem 1.4]{Tan1}. 
On the other hand, even if $k$ is an imperfect field, 
the $p$-power freeness of $L$ holds. 
Indeed, if either $\dim Y \geq 1$ or $\dim Y=0$ and $L \not\equiv 0$, 
then the assertion follows from \cite[Theorem 4.2]{Tan18}. 
If $\dim Y=0$ and $L \equiv 0$, then 
we obtain $p^e L \sim 0$ by \cite[Theorem 1.3]{BT}. 
We do not use this fact in this paper. 
\end{rem}

\begin{thm}\label{t-birat-torsion}
Let $k$ be a perfect field of characteristic $p>5$. 
Let $f:X \to Y$ be a projective birational $k$-morphism of quasi-projective normal threefolds over $k$. 
Assume that there exists an effective $\R$-divisor $\Delta_Y$ on $Y$ such that 
$(Y, \Delta_Y)$ is klt. 
Let $L$ be a Cartier divisor on $X$ such that $L \equiv_f 0$. 
Then there exist a positive integer $e$ and a Cartier divisor $L_Y$ on $Y$ such that $L \sim p^eL_Y$. 
\end{thm}

\begin{proof}
%It is enough to show that $L$ is $p$-power free over $Y$. 
%Thus, we may assume that $Y$ is affine. In particular, $X$ and $Y$ are quasi-projective over $k$. 
The proof consists of three steps. 

\setcounter{step}{0}

\begin{step}\label{s1-birat-torsion}
The assertion of Theorem \ref{t-birat-torsion} holds if $Y$ is $\Q$-factorial.
\end{step}

\begin{proof}[Proof of Step \ref{s1-birat-torsion}]
Taking a log resolution of $(Y, \Delta_Y)$ which dominates $X$, 
we may assume that $f:X \to Y$ is a log resolution of $(Y, \Delta_Y)$. 
Let $E$ be the reduced $f$-exceptional divisor such that $\Supp\,E=\Ex(f)$. 
Set $\Delta:=f_*^{-1}\Delta_Y + E$. 
Then we have 
\[
K_X+\Delta=K_X+f_*^{-1}\Delta_Y + E \sim_{\R} f^*(K_Y+\Delta_Y)+ E'
\] 
where $E'$ is an $f$-exceptional effective $\R$-divisor on $X$ such that 
$\Supp\,E'=\Ex(f)$. 
By \cite[Theorem 1.1]{HNT}, 
there exists a $(K_X+\Delta)$-MMP over $Y$ that terminates: 
\[
X =: X_0 \dashrightarrow X_1 \dashrightarrow \cdots \dashrightarrow X_N. 
\]
The negativity lemma implies that $f_N:X_N \to Y$ is small. 
Since $Y$ is $\Q$-factorial, $f_N$ is an isomorphism. 

Therefore, it is enough to treat the case when 
$(X, \Delta)$ is a $\Q$-factorial dlt pair, 
$-(K_X+\Delta)$ is $f$-ample, $\rho(X/Y)=1$, and $-S$ is $f$-ample for some prime divisor $S$ contained in $\Supp\,\llcorner \Delta \lrcorner$. 
It holds that $A:=L-S$ is an $f$-ample $\Q$-Cartier $\Z$-divisor. 
Since $L|_S$ is $p$-power free over $Y$ (Lemma \ref{l-pp-2dim}), 
it follows from Proposition \ref{p-Keel-div} that $L$ is $p$-power free over $Y$. 
This completes the proof of Step \ref{s1-birat-torsion}. 
\end{proof}

\begin{step}\label{s2-birat-torsion}
The assertion of Theorem \ref{t-birat-torsion} holds 
if there exists an effective $\Q$-divisor $\Delta$ on $X$ such that 
$(X, \Delta)$ is klt and $-(K_X+\Delta)$ is $f$-nef and $f$-big. 
\end{step}

\begin{proof}[Proof of Step \ref{s2-birat-torsion}]
Set $\mcL:=\MO_X(L)$. 
Taking the base change to an uncountable algebraically closed field, 
we may assume that $k$ is an uncountable algebraically closed field 
(Lemma \ref{l-pp-ff-descent}). 
For an arbitrary closed point $y \in Y$, 
it is enough to show that $\mcL|_{X_y}$ is $p$-power free 
(Theorem \ref{t-rel-pp-free2}). 
By \cite[Proposition 2.15]{GNT}, 
there exists a commutative diagram consisting of projective birational morphisms 
of quasi-projective normal threefolds 
$$\begin{CD}
W @>\psi >> X'\\
@VV\varphi V @VVf'V\\
X @>f>> Y
\end{CD}$$
and an effective $\Q$-divisor $\Delta'$ on $X'$ which satisfy the following properties. 
\begin{enumerate}
\item $(X', \Delta')$ is a $\Q$-factorial plt threefold and $-(K_{X'}+\Delta')$ is $f'$-ample. 
\item $-\llcorner \Delta' \lrcorner$ is $f'$-nef and $f'^{-1}(z)_{\red} = \llcorner \Delta'\lrcorner$. 
%\item 
%$W$ is a smooth threefold, and both of $\varphi$ and $\psi$ are projective birational morphisms. 
\item 
For $T:=\llcorner \Delta' \lrcorner$, 
there exists an effective $\Q$-divisor $\Delta_T$ on $T$ such that 
$(T, \Delta_T)$ is klt and $-(K_T+\Delta_T)$ is ample. 
\end{enumerate}
Set $\mcL_W:=\varphi^*\mcL$. 
After replacing $\mcL$ by $\mcL^{\otimes p^e}$ for some $e \in \Z_{>0}$, 
Step \ref{s1-birat-torsion} enables us to find 
an invertible sheaf $\mcL'$ on $X'$ such that $\mcL_W \simeq \psi^*\mcL'$. 
Then we have that $\mcL'|_{(X'_y)_{\red}}=\mcL'|_{T}$ 
is $p$-power free (Lemma \ref{l-pp-2dim}). 
It follows from Lemma \ref{l-pp-univ-homeo} that 
also $\mcL'|_{X'_y}$ is $p$-power free. 
Since $\mcL_W|_{W_y}$ is the pullback of $\mcL'|_{X'_y}$, 
it holds that $\mcL_W|_{W_y}$ is $p$-power free. 
Finally, since $W_y \to X_y$ has connected fibres, 
it follows again from Lemma \ref{l-pp-univ-homeo} that 
$\mcL|_{X_y}$ is $p$-power free. 
This completes the proof of Step \ref{s2-birat-torsion}. 
\end{proof}

\begin{step}\label{s3-birat-torsion}
The assertion of Theorem \ref{t-birat-torsion} holds without any additional assumptions.
\end{step}

\begin{proof}[Proof of Step \ref{s3-birat-torsion}]
By the same argument as in  Step \ref{s1-birat-torsion}, 
the problem is reduced to the case when 
$X$ is $\Q$-factorial and $f$ is a small birational morphism. 
In particular, for $\Delta:=f_*^{-1}\Delta_Y$, 
it holds that $K_X+\Delta = f^*(K_Y+\Delta_Y)$. 
In particular, $-(K_X+\Delta)$ is $f$-nef and $f$-big. 
%Then there exists an effective $\Q$-divisor $\Delta'$ on $X$ such that $(X, \Delta')$ is klt and $-(K_X+\Delta')$ is $f$-ample. 
Then we may apply Step \ref{s2-birat-torsion}. 
This completes the proof of Step \ref{s3-birat-torsion}. 
\end{proof}
Step \ref{s3-birat-torsion} completes the proof of Theorem \ref{t-birat-torsion}. 
\end{proof}

\begin{thm}\label{t-bpf-pp}
Let $k$ be a perfect field of characteristic $p>5$. 
Let $(X, \Delta)$ be a three-dimensional klt pair and 
let $f:X \to Z$ be a projective surjective $k$-morphism to a quasi-projective $k$-scheme $Z$. 
Let $L$ be an $f$-nef Cartier divisor on $X$ such that $L-(K_X+\Delta)$ is $f$-nef and $f$-big. 
Assume that either 
\begin{enumerate}
\item $\dim Z \geq 1$, or 
\item $\dim Z=0$ and $L \not\equiv 0$. 
\end{enumerate}
Then $L$ is $p$-power $f$-free.  
\end{thm}

\begin{proof}
The proof consists of two steps. 

\setcounter{step}{0}
\begin{step}\label{s1-bpf-pp}
The assertion of Theorem \ref{t-bpf-pp} holds if (1) holds. 
\end{step}

\begin{proof}[Proof of Step \ref{s1-bpf-pp}] 
Set $\mcL:=\MO_X(L)$. 
Then, by Lemma \ref{l-pp-ff-descent}, 
we may assume that $k$ is an uncountable algebraically closed field. 

Fix a closed point $z \in Z$. 
By \cite[Proposition 2.15]{GNT}, 
there exists a commutative diagram consisting of projective morphisms 
of quasi-projective normal varieties 
$$\begin{CD}
W @>\psi >> X'\\
@VV\varphi V @VVf'V\\
X @>f>> Z
\end{CD}$$
and an effective $\Q$-divisor $\Delta'$ on $X'$ which satisfy the following properties. 
\begin{enumerate}
\item $(X', \Delta')$ is a $\Q$-factorial plt threefold and $-(K_{X'}+\Delta')$ is $f'$-ample. 
\item $-\llcorner \Delta' \lrcorner$ is $f'$-nef and $f'^{-1}(z)_{\red} = \llcorner \Delta'\lrcorner$. 
\item 
$W$ is a smooth threefold, and both of $\varphi$ and $\psi$ are projective birational morphisms. 
\item 
For $T:=\llcorner \Delta' \lrcorner$, 
there exists an effective $\Q$-divisor $\Delta_T$ on $T$ such that 
$(T, \Delta_T)$ is klt and $-(K_T+\Delta_T)$ is ample. 
\end{enumerate}
Then, by the same argument as in Step \ref{s2-birat-torsion} of 
Theorem \ref{t-birat-torsion}, it holds that $\mcL|_{X_z}$ 
is $p$-power free for any closed point $z \in Z$. 
Then it follows from Theorem \ref{t-rel-pp-free2} that 
$\mcL$ is $p$-power $f$-free. 
This completes the proof of Step \ref{s1-bpf-pp}
\end{proof}

\begin{step}\label{s2-bpf-pp}
The assertion of Theorem \ref{t-bpf-pp} holds if (2) holds. 
\end{step}

\begin{proof}[Proof of Step \ref{s2-bpf-pp}] 
It follows from \cite[Theorem 2.9]{GNT} that $L$ is $f$-semi-ample. 
Hence, there exist a positive integer $m$, 
projective morphisms 
\[
f:X \xrightarrow{g} Y \xrightarrow{f_Y} Z, 
\] 
and an ample Cartier divisor $A$ on $Y$ such that 
$g_*\MO_X=\MO_Y$ and $mL \sim g^*A$. 
In particular, we have $L \equiv_g 0$. 
By (2), it holds that $\dim Y >0$. 
By Step \ref{s1-bpf-pp}, 
there exist $e \in \Z_{>0}$ and an $f_Y$-ample Cartier divisor $A'$ on $Y$ 
such that $p^eL \sim g^*A'$. 
This completes the proof of Step \ref{s2-bpf-pp}. 
\end{proof}
Step \ref{s1-bpf-pp} and Step \ref{s2-bpf-pp} complete the proof of Theorem \ref{t-bpf-pp}. 
\end{proof}

%\footnote{aaa}

%%%%%%%%%%%%%%%%%%%%%%%%%%%%%%%%%%%%%%%%%%%%%%%%%%%%%%%%%%%%%%%%%%%%

\end{document}